\documentclass[12pt]{amsart}

\numberwithin{equation}{section}

\newtheorem{theorem}[equation]{Theorem}
\newtheorem{lemma}[equation]{Lemma}
\newtheorem{corollary}[equation]{Corollary}
\newtheorem{proposition}[equation]{Proposition}

\newtheorem{mainconj}[equation]{The Main Conjecture}

\theoremstyle{definition}
\newtheorem{definition}[equation]{Definition}

\theoremstyle{remark}

\newenvironment{pf}
{\noindent\textbf{Proof: }}
{\hfill\textbf{ Q.E.D.}\medskip}

\def\LS#1#2#3#4{\mathcal{L}_{#1}(#2,{#4}^{#3})}
\def\LSH#1#2#3{\mathcal{L}_{#1}({#3}^{#2})}

\newcommand{\bP}{{\mathbb{P}}}
\newcommand{\bF}{{\mathbb{F}}}
\renewcommand{\L}{\mathcal{L}}
\newcommand{\M}{\mathcal{M}}
\newcommand{\Lo}{\mathcal{L}_0}

\newcommand{\LP}{\mathcal{L}_\bP}
\newcommand{\LF}{\mathcal{L}_\bF}
\newcommand{\LPhat}{\hat{\mathcal{L}}_\mathbb{P}}
\newcommand{\LFhat}{\hat{\mathcal{L}}_\mathbb{F}}

\renewcommand{\l}{\ell}
\newcommand{\lo}{\ell_0}

\newcommand{\lP}{\ell_\mathbb{P}}
\newcommand{\lF}{\ell_\mathbb{F}}
\newcommand{\lPhat}{\hat{\ell}_\mathbb{P}}
\newcommand{\lFhat}{\hat{\ell}_\mathbb{F}}
\newcommand{\rP}{r_\mathbb{P}}
\newcommand{\rF}{r_\mathbb{F}}
\renewcommand{\v}{v}

\newcommand{\vP}{v_\mathbb{P}}
\newcommand{\vF}{v_\mathbb{F}}
\newcommand{\vPhat}{\hat{v}_\mathbb{P}}
\newcommand{\vFhat}{\hat{v}_\mathbb{F}}
\newcommand{\e}{e}
\renewcommand{\O}{\mathcal{O}}
\newcommand{\Pic}{\mathrm{Pic}}
\newcommand{\dhigh}{\operatorname{dhigh}}
\newcommand{\dlow}{\operatorname{dlow}}

\title{Linear Systems of Plane Curves
with Base Points of Equal Multiplicity}
\author{Ciro Ciliberto and Rick Miranda}
\address{Dept. of Mathematics\\Universit\'a di Roma II\\Via Fontanile di
Carcaricola\\00173 Rome, Italy}
\email{cilibert@axp.mat.utovrm.it}
\address{Dept. of Mathematics\\Colorado State University\\Ft. Collins, CO
80523}
\email{miranda@math.colostate.edu}
\thanks{Research supported in part by the NSA}
\thanks{\today. Set in type by LaTeX2e and amstex}

\begin{document}

\begin{abstract}
In this article we address the problem of computing
the dimension of the space of plane curves of degree $d$
with $n$ general points of multiplicity $m$.
A conjecture of Harbourne \cite{harbourne1}
and Hirschowitz \cite{hirschowitz3}
implies that when $d \geq 3m$,
the dimension is equal to the expected dimension
given by the Riemann-Roch Theorem.
Also, systems for which the dimension is larger than expected
should have a fixed part containing a multiple $(-1)$-curve.
We reformulate this conjecture by explicitly listing
those systems which have unexpected dimension.
Then we use a degeneration technique
developed in \cite{cm1}
to show that the conjecture holds
for all $m \leq 12$.
\end{abstract}

\maketitle

\tableofcontents

\section*{Introduction}
Consider the projective plane $\bP^2$
and $n+1$ general points $p_0, p_1,\ldots, p_n$ on it.
Let $H$ denote the line class of the plane.
Consider the linear system consisting of
plane curves of degree $d$ (that is, divisors in $|dH|$)
with multiplicity at least $m_i$ at $p_i$ for $i \geq 0$.
If all $m_i$ for $i \geq 1$ are equal, to $m$ say,
we denote this system by $\L = \LS{d}{m_0}{n}{m}$
and call the system \emph{quasi-homogeneous}.
If in addition $m_0=0$, we say the system is \emph{homogeneous},
and denote it simply by $\LSH{d}{n}{m}$.

Define its \emph{virtual dimension}
\[
\v = d(d+3)/2 - m_0(m_0+1)/2 - nm(m+1)/2;
\]
the actual dimension of the linear system
cannot be less than $-1$,
and hence we define the \emph{expected dimension} to be
\[
\e = \max\{-1, \v \}.
\]

The dimension of $\L$ achieves its minimum value
for a general set of points;
abusing notation slightly we call this the \emph{dimension} of $\L$,
and denote it by $\l$.
We always have that
\[
\l \geq \e;
\]
equality implies
(when the numbers are at least $-1$)
that the conditions imposed by the multiple points are independent.

We will say that the system $\L$ is {\em non-special}
if equality holds, i.e.,
that either the system is empty
or that the conditions imposed by the multiple points are independent.
If $\l > \e$ then we say the system is {\em special}.

The speciality of $\L$
is equivalent to a statement about linear systems
on the blowup $\bP'$ of $\bP^2$ at the points $p_i$.
If $H$ denotes the class of the pullback of a line
and $E_i$ denotes the class of the exceptional divisor above $p_i$,
then the linear system on $\bP^2$
transforms to the linear system
$\L' = |dH-\sum_{i=0}^n m_iE_i|$ on $\bP'$.
Then the original system $\L$ is non-special if and only if
\[
h^1(\L') = \max\{0, -1-\v \}.
\]
In particular if the system is non-empty
(which means that $H^0(\L')$ is non-zero),
it is non-special if and only if the $H^1$ is zero.
More precisely if the virtual dimension $\v \geq -1$
then non-speciality means that the $H^1$ is zero,
or, equivalently, that the conditions imposed by
the multiple base points are linearly independent.

The \emph{self-intersection} ${\L}^2$ and the
\emph{genus} $g_{\L}$
are defined in terms of the transformed system $\L'$
on the blowup;
we have:
\[
{\L}^2=d^2-m^2_0-nm^2
\;\;\;\mathrm{ and }\;\;\;
2g_{\L}-2=d(d-3)-m_0(m_0-1)-nm(m-1).
\]
Notice the basic identity:
\begin{equation}
\label{RR}
{\v}={\L}^2-g_{\L}+1.
\end{equation}

The \emph{intersection number}
${\L}(d,m_0,n,m)\cdot{\L}(d',m_0',n',m')$, $n'\leq n$,
is given by:
\[
{\L}(d,m_0,n,m)\cdot{\L}(d',m'_0,n',m) := dd' - m_0m_0' - n'mm'
\]

\section{The Degeneration of the Plane and the Recursion}

In this section we describe the degeneration of the plane
which we use in the analysis;
for details of the construction we refer to \cite{cm1}.

Let $\Delta$ be a complex disc around the origin.
The product $V = \bP^2 \times \Delta$
comes with the two projections
$p_1:V \to \Delta$ and $p_2:V \to \bP^2$;
let $V_t = \bP^2\times \{t\}$.

Blow up a line $L$ in the plane $V_0$
and obtain a new three-fold $X$ with maps
$f:X \to V$,
$\pi_1 = p_1 \circ f: X \to \Delta$,
and $\pi_2 = p_2 \circ f: X \to \bP^2$.
The map $\pi_1:X \to \Delta$
is a flat family of surfaces $X_t = \pi_1^{-1}(t)$ over $\Delta$.
If $t \neq 0$, then $X_t = V_t$ is a plane $\bP^2$,
while $X_0$ is the union
of the proper transform $\bP$ of $V_0$
and of the exceptional divisor $\bF$ of the blow-up.
The surface $\bP$ is a plane $\bP^2$
and $\bF$ is a Hirzebruch surface $\bF_1$.
They are joined transversally along a curve $R$
which is a line $L$ in $\bP$
and is the exceptional divisor $E$ on $\bF$.

The Picard group of $X_0$
is the fibered product of $\Pic(\bP)$ and $\Pic(\bF)$
over $\Pic(R)$;
a line bundle $\mathcal{X}$ on $X_0$ is a line bundle
$\mathcal{X}_{\bP}$ on $\bP$
and a line bundle $\mathcal{X}_{\bF}$ on $\bF$
which agree on the double curve $R$.
This means that
$\mathcal{X}_{\bP}\cong{\O}_{\bP}(d)$ and
$\mathcal{X}_{\bF}\cong{\O}_{\bF}(cH-dE)$ for some $c$ and $d$.
We will denote this line bundle by $\mathcal{X}(c,d)$.

Note that the bundle $\O_X(\bP)$
restricts to $\bP$ as ${\O}_{\bP}(-1)$
and restricts to $\bF$ as ${\O}_{\bF}(E)$.
Let ${\O}_X(d)$ be the line bundle $\pi_2^*({\O}_{\bP^2}(d))$;
it restricts to ${\O}_{\bP}(d)$ and to ${\O}_{\bF}(dH-dE)$ on $\bF$.
Let us denote by ${\O}_X(d,a)$ the line bundle
${\O}_X(d)\otimes {\O}_X((d-a)\bP)$.
The restriction of ${\O}_X(d,a)$ to $X_t$, $t\not=0$,
is isomorphic to ${\O}_{\bP^2}(d)$,
but the restriction to $X_0$ is 
isomorphic to $\mathcal{X}(d,a)$
We therefore see that all of the bundles $\mathcal{X}(d,a)$ on $X_0$
are flat limits of the bundles $\O_{\bP^2}(d)$
on the general fiber $X_t$ of this degeneration.

Fix a positive integer $n$
and another non-negative integer $b\leq n$.
Let us consider $n-b+1$ general points $p_0, p_1, \ldots, p_{n-b}$ in $\bP$
and $b$ general points $p_{n-b+1},...,p_n$ in $\bF$.
These points are limits
of $n$ general points $p_{0,t}, p_{1,t}, \ldots, p_{n,t}$ in $X_t$.
Consider then the linear system  ${\L}_t$
which is the system $\LS{d}{m_0}{n}{m}$ in $X_t \cong {\bP}^2$
based at the points $p_{0,t}, p_{1,t},...,p_{n,t}$.

We now also consider the linear system
${\Lo}$ on $X_0$
which is formed by the divisors in $|\mathcal{X}(d,a)|$
having a point of multiplicity $m_0$ at $p_0$
and points of multiplicity $m$ at $p_1,...,p_n$.
According to the above considerations,
any one of the systems ${\Lo}$ (for any $a$ and $b$)
can be considered as a flat limit on $X_0$
of the system ${\L}=\LS{d}{m_0}{n}{m}$.
We will say that ${\Lo}$ is obtained from $\L$
by an \emph{$(a,b)$-degeneration}.

We note that the system $\Lo$ restricts to $\bP$
as a system $\LP$ of the form $\LS{a}{m_0}{n-b}{m}$
and $\Lo$ restricts to $\bF$
as a system $\LF$ of the form $\LS{d}{a}{b}{m}$.
Indeed, at the level of vector spaces,
the system $\Lo$ is the fibered product
of $\LP$ and $\LF$ over the restricted system on $R$,
which is $\O_R(a)$.

We denote by $\lo$ the dimension of the linear system $\Lo$ on $X_0$.
By semicontinuity, this dimension $\lo$
is at least that of the linear system on the general fiber,
i.e.,
\[
\lo = \dim(\Lo) \geq \l = \dim \LS{d}{m_0}{n}{m}.
\]
Therefore we have the following:

\begin{lemma}
\label{lo=Ethenreg}
If $\lo$ is equal to the expected dimension $\e$
of $\L=\LS{d}{m_0}{n}{m}$
then the system $\L$ is non-special.
\end{lemma}

The main result of \cite{cm1} was the computation of the
dimension $\lo$ of the limit linear system $\Lo$.
We will not reproduce the argument here.
The dimension $\lo$ is obtained in terms of the
dimensions of the systems $\LP$ and $\LF$,
and the dimensions of the subsystems $\LPhat \subset \LP$
and $\LFhat \subset \LF$
consisting of divisors containing the double curve $R$.
Notice that by slightly abusing notation we have
\begin{align*}
\LP &= \LS{a}{m_0}{n-b}{m} \\
\LPhat &= \LS{a-1}{m_0}{n-b}{m} \\
\LF &= \LS{d}{a}{b}{m}, \\
\LFhat &= \LS{d}{a+1}{b}{m};
\end{align*}
all of these systems are quasi-homogeneous,
which provides the basis for the recursion.

We denote by
\[
\begin{array}{ll}
\vP,\;\;\vF &
\text{ the virtual dimension of the systems} \;\;\LP,\;\; \LF \\
\vPhat,\;\;\vFhat &
\text{ the virtual dimension of the subsystems} \;\;\LPhat, \;\;\LFhat \\
\lP,\;\;\lF &
\text{ the dimension of the systems} \;\;\LP, \;\;\LF \\
\lPhat,\;\;\lFhat&
\text{ the dimension of the subsystems} \;\;\LPhat, \;\;\LFhat. \\
\end{array}
\]
The following lemma gives three useful identities
(of polynomials in $d$, $m_0$, $n$, $m$, $a$, and $b$)
which the reader can easily check.

\begin{lemma}
\label{vP+vF}
$\v = \vP+\vF - a = \vF + \vPhat +1 = \vP + \vFhat +1$.
\end{lemma}

\begin{theorem}
\label{dimLo}
Let $\rP = \lP-\lPhat-1$ and $\rF = \lF-\lFhat-1$;
these are the dimensions of the restrictions (to $R$)
of the linear systems $\LP$ and $\LF$ respectively.
\newline(a) If $\rP+\rF \leq a-1$, then
$\lo = \lPhat + \lFhat + 1$.
\newline(b) If $\rP+\rF \geq a-1$, then
$\lo = \lP+ \lF - a$.
\end{theorem}

\section{Homogeneous $(-1)$-Configurations%
and the Main Conjecture}

A linear system $\L = \LS{d}{m_0}{n}{m}$
with $\L^2 = -1$ and $g_\L = 0$
will be called a \emph{quasi-homogeneous $(-1)$-class}.
By (\ref{RR}), we see that $\v = 0$,
so that every quasi-homogeneous $(-1)$-class is effective.

Suppose that $A$ is an irreducible rational curve
and is a member of a linear system $\L=\LS{d}{m_0}{n}{m}$.
If on the blowup $\bP^\prime$ of the plane
the proper transform of $A$ is smooth, of self-intersection $-1$,
then we say that $A$ is a \emph{$(-1)$-curve}.
Such a linear system $\L$ is non-special, of dimension $0$.
A quasi-homogeneous $(-1)$-class
containing a $(-1)$-curve
will be called an \emph{irreducible $(-1)$-class}.

\begin{definition}
A linear system $\L$ is \emph{$(-1)$-special}
if there are $(-1)$-curves $A_1,\ldots,A_r$
such that $\L\cdot A_j = -N_j$ with $N_j \geq 1$ for every $j$
and $N_j \geq 2$ for some $j$,
with the residual system $\M = \L - \sum_j N_j A_j$
having non-negative virtual dimension $\v(\M) \geq 0$,
and having non-negative intersection with every $(-1)$-curve.
\end{definition}

We refer the reader to \cite{cm1} for comments on this definition.
We note in particular
that if there are different $(-1)$-curves $A_i$ and $A_j$
both of which meet $\L$ negatively,
then it must be the case that $A_i \cdot A_j = 0$.
A divisor which is the sum of pairwise disjoint $(-1)$-curves
will be called a \emph{$(-1)$-configuration}.

Every $(-1)$-special system is special;
see Lemma 4.1 of \cite{cm1}.
The main conjecture that we are concerned with
is the following restatement of a conjecture
of Hirschowitz (see \cite{hirschowitz3}).

\begin{mainconj}
Every special system is $(-1)$-special.
\end{mainconj}

Suppose a quasi-homogeneous system $\LS{d}{m_0}{n}{m}$
meets negatively a $(-1)$-curve $A$
of degree $\delta$, having multiplicities
$\mu_0,\mu_1,\ldots,\mu_n$ at the points $p_0,\ldots,p_n$.
Since the points are general, by monodromy
we have that for any permutation $\sigma \in \Sigma_n$,
$\L$ also meets negatively
the $(-1)$-curve $A_\sigma$ of degree $\delta$,
having multiplicity $\mu_0$ at $p_0$,
and having multiplicities $\mu_{\sigma(i)}$ at $p_i$
for each $i \geq 1$.
Thus either $A$ is quasi-homogeneous itself
(and all $A_\sigma$'s are equal to $A$)
or we obtain a $(-1)$-configuration formed by $A$
and the other $A_\sigma$'s.
Necessarily, this configuration is quasi-homogeneous,
and is homogeneous if the original linear system $\L$ is.

In \cite{cm1}, Proposition 5.16,
we classified all homogeneous $(-1)$-configurations;
this list we reproduce below; see also Theorem 7 of \cite{nagata}.

\begin{proposition}
\label{homog-1}
The following is a complete list of homogeneous $(-1)$-configurations:
\begin{align*}
\LSH{1}{2}{1}: &\;\; \text{ a line through $2$ points} \\
\LSH{2}{5}{1}: &\;\; \text{ a conic through $5$ points}\\
\LSH{3}{3}{2}: &\;\; \text{ $3$ lines each passing through $2$ of $3$ points} \\
\LSH{12}{6}{5}: &\;\; \text{ $6$ conics each passing through $5$ of $6$ points}\\
\LSH{21}{7}{8}: &\;\; \text{ $7$ cubics through $6$ points, double at another}\\
\LSH{48}{8}{17}: &\;\; \text{ $8$ sextics double at $7$ points, triple at another} \\
\end{align*}
\end{proposition}

This list enables us to classify all homogeneous $(-1)$-special systems;
they are the systems $\LSH{d}{n}{m}$ which intersect
one of the above curves negatively.

\begin{theorem}
\label{list-1special}
The $(-1)$-special homogeneous linear systems are
\begin{align*}
\LSH{d}{2}{m} &\;\;\text{with}\;\; m      \leq d \leq 2m-2, \\
\LSH{d}{3}{m} &\;\;\text{with}\;\; 3m/2   \leq d \leq 2m-2, \\
\LSH{d}{5}{m} &\;\;\text{with}\;\; 2m     \leq d \leq (5m-2)/2, \\
\LSH{d}{6}{m} &\;\;\text{with}\;\; 12m/5  \leq d \leq (5m-2)/2, \\
\LSH{d}{7}{m} &\;\;\text{with}\;\; 21m/8  \leq d \leq (8m-2)/3,\\
\LSH{d}{8}{m} &\;\;\text{with}\;\; 48m/17 \leq d \leq (17m-2)/6.\\
\end{align*}
\end{theorem}

\begin{pf}
Let $\L=\LSH{d}{n}{m}$ be a homogeneous $(-1)$-special system.
Then $\L=\M+NA$, where $A$ is a homogeneous $(-1)$-configuration
and $\M$ is a homogeneous system
with $\v(\M) \geq 0$ and $\M\cdot A = 0$.
We have exactly six possibilities for $A$,
and hence for $n$,
given Proposition \ref{homog-1}.
We take these up in turn,
seeking the homogeneous system $\M$.

\medskip
\noindent{$\underline{A = \LSH{1}{2}{1}}$:}
Let $\M=\LSH{\delta}{2}{\mu}$.
The condition that $\M\cdot A = 0$ is that $\delta=2\mu$,
so that $\M=\LSH{2\mu}{2}{\mu}$;
then $\v(\M) = \mu(\mu+2)$ which is always $\geq 0$ if $\mu$ is.
Therefore the $(-1)$-special system $\L$ is of the form
$\L = \M+NA=\LSH{2\mu+N}{2}{(\mu+N)}$ with $N \geq 2$ and $\mu\geq 0$.
This is a general homogeneous system $\LSH{d}{2}{m}$
with $m \leq d \leq 2m-2$.

\medskip
\noindent{$\underline{A=\LSH{2}{5}{1}}$:}
Let $\M=\LSH{\delta}{5}{\mu}$.
The condition that $\M\cdot A = 0$ is that $2\delta=5\mu$,
so that there is an integer $k\geq 0$ with $\mu=2k$ and $\delta=5k$,
and therefore $\M=\LSH{5k}{5}{(2k)}$.
Then $\v(\M) = 5k(k+1)/2$ which is always $\geq 0$ if $k$ is.
Therefore the $(-1)$-special system $\L$ is of the form
$\L = \M+NA=\LSH{5k+2N}{5}{(2k+N)}$ with $N \geq 2$ and $k\geq 0$.
This is a general homogeneous system $\LSH{d}{5}{m}$
with $2m \leq d \leq (5m-2)/2$.

\medskip
\noindent{$\underline{A=\LSH{3}{3}{2}}$:}
Here $A$ consists of three curves,
and the $(-1)$-curve is $A_0\in \LSH{1}{2}{1}$,
a line through $2$ of the three points.
Let $\M=\LSH{\delta}{3}{\mu}$.
The condition that $\M\cdot A = 0$ is that $\delta=2\mu$,
so that $\M=\LSH{2\mu}{3}{\mu}$;
then $\v(\M) = \mu(\mu+3)/2$ which is always $\geq 0$ if $\mu$ is.
Therefore the $(-1)$-special system $\L$ is of the form
$\L = \M+NA=\LSH{2\mu+3N}{3}{(\mu+2N)}$ with $N \geq 2$ and $\mu\geq 0$.
This is a general homogeneous system $\LSH{d}{3}{m}$
with $3m/2 \leq d \leq 2m-2$.

\medskip
\noindent{$\underline{A=\LSH{12}{6}{5}}$:}
Here $A$ consists of $6$ conics,
and the $(-1)$-curve is $A_0\in \LSH{2}{5}{1}$,
a conic through $5$ of the $6$ points.
Let $\M=\LSH{\delta}{6}{\mu}$.
The condition that $\M\cdot A = 0$ is that $2\delta=5\mu$,
so that there is an integer $k\geq 0$ with $\mu=2k$ and $\delta=5k$,
and therefore $\M=\LSH{5k}{6}{(2k)}$.
Then $\v(\M) = k(k+3)/2$ which is always $\geq 0$ if $k$ is.
Therefore the $(-1)$-special system $\L$ is of the form
$\L = \M+NA=\LSH{5k+12N}{6}{(2k+5N)}$ with $N \geq 2$ and $k\geq 0$.
This is a general homogeneous system $\LSH{d}{6}{m}$
with $12m/5 \leq d \leq (5m-2)/2$.

\medskip
\noindent{$\underline{A=\LSH{21}{7}{8}}$:}
Here $A$ consists of $7$ cubics,
and the $(-1)$-curve is $A_0\in \LS{3}{2}{6}{1}$,
a cubic double at one point
and passing through the other $6$ of the $7$ points.
Let $\M=\LSH{\delta}{7}{\mu}$;
the condition that $\M\cdot A = 0$ is that $3\delta=8\mu$,
so that there is an integer $k\geq 0$ with $\mu=3k$ and $\delta=8k$,
and therefore $\M=\LSH{8k}{7}{(3k)}$.
Then $\v(\M) = k(k+3)/2$ which is always $\geq 0$ if $k$ is.
Therefore the $(-1)$-special system $\L$ is of the form
$\L = \M+NA=\LSH{8k+21N}{7}{(3k+8N)}$ with $N \geq 2$ and $k\geq 0$.
This is a general homogeneous system $\LSH{d}{7}{m}$
with $21m/8 \leq d \leq (8m-2)/3$.

\medskip
\noindent{$\underline{A=\LSH{48}{8}{17}}$:}
Here $A$ consists of $8$ sextics,
and the $(-1)$-curve is $A_0\in \LS{6}{3}{7}{2}$,
a sextic triple at one point
and double at the other $7$ of the $8$ points.
Let $\M=\LSH{\delta}{8}{(m')}$;
the condition that $\M\cdot A = 0$ is that $6\delta=17m'$,
so that there is an integer $k\geq 0$ with $m'=6k$ and $\delta=17k$,
and therefore $\M=\LSH{17k}{8}{(6k)}$.
Then $\v(\M) = k(k+3)/2$ which is always $\geq 0$ if $k$ is.
Therefore the $(-1)$-special system $\L$ is of the form
$\L = \M+NA=\LSH{17k+48N}{8}{(6k+17N)}$ with $N \geq 2$ and $k\geq 0$.
This is a general homogeneous system $\LSH{d}{8}{m}$
with $48m/17 \leq d \leq (17m-2)/6$.
\end{pf}

The reader will notice that the only $(-1)$-special homogeneous
systems occur when $n \leq 8$.  In fact the Main Conjecture is true
in this range, which is a classical fact, see
\cite{nagata}, \cite{harbourne1}, \cite{gimigliano1}, \cite{gimigliano2}.
The precise statement we will find useful is the following.

\begin{theorem}
The Main Conjecture is true for all homogeneous linear systems
$\LSH{d}{n}{m}$ with $n \leq 9$.
In particular every homogeneous linear system $\LSH{d}{4}{m}$
and $\LSH{d}{9}{m}$ is non-special.
\end{theorem}

\section{The Recursion for Homogeneous Systems}
Suppose that we want to investigate the dimension 
of a homogeneous system $\L=\LSH{d}{n}{m}$.
We construct an $(a,b)$ degeneration of the plane and the bundle,
and we are led to studying the four systems
\begin{align*}
\LP &= \LS{a}{0}{n-b}{m} \\
\LPhat &= \LS{a-1}{0}{n-b}{m} \\
\LF &= \LS{d}{a}{b}{m}, \\
\LFhat &= \LS{d}{a+1}{b}{m}.
\end{align*}
We note that the first two systems $\LP$ and $\LPhat$
are also homogeneous, and so an opportunity to apply
induction presents itself.
However the last two systems $\LF$ and $\LFhat$ on $\bF$
are not homogeneous in general,
which spoils the possibility of a simple-minded induction
on this side.

The reader sees that we need a second method
to compute the dimensions of these last two systems.
Such a method is provided by using Cremona transformations
when the extra multiplicity ($a$ and $a+1$ in the case
of $\LF$ and $\LFhat$ respectively)
are large with respect to the degree.
In \cite{cm1} we have made the analysis necessary
and we present the results below.

We first handle the case of a linear system
$\LS{d}{a}{b}{m}$ with $a = d-m$.

\begin{proposition}
\label{m0=d-m}
Let $\L=\LS{d}{d-m}{b}{m}$ with $2 \leq m \leq d$.
Write $d = qm+\mu$ with $0\leq \mu\leq m-1$,
and $b = 2h+\epsilon$, with $\epsilon \in\{0,1\}$.
Then the system $\L$ is special if and only if
$q = h$, $\epsilon = 0$, and $\mu \leq m-2$.
More precisely:
\begin{itemize}
\item[(a)] If $q \geq h+1$ then $\L$ is nonempty and non-special.
In this case 
\[
\dim\L=d(m+1)-\binom{m}{2}-b\binom{m+1}{2}.
\]
\item[(b)] If $q=h$ and $\epsilon = 1$
the system $\L$ is empty and non-special.
\item[(c)] If $q=h$, $\epsilon = 0$, and $\mu = m-1$,
the system $\L$ is nonempty and non-special;
in this case
\[
\dim\L=(m-1)(m+2)/2.
\]
\item[(d)] If $q=h$, $\epsilon = 0$, and $\mu \leq m-2$,
the system $\L$ is special;
in this case
\[
\dim \L = \mu(\mu+3)/2.
\]
\item[(e)] If $q\leq h-1$
the system $\L$ is empty and non-special.
\end{itemize}
\end{proposition}

If $a > d-m$, then in the system $\LS{d}{a}{b}{m}$,
the lines through $p_0$ and $p_i$ split off
repeatedly, with a residual system having $m_0 = d-m$.
Therefore the above analysis leads to a computation
in these cases also.

\begin{corollary}
\label{m0=d-m+k}
Let $\L = \LS{d}{d-m+k}{b}{m}$ with $k \geq 1$,
and let
\[
\L' = \LS{d-kb}{d-kb-m+k}{b}{m-k}.
\]
Then $\dim \L = \dim \L'$ and
$\L$ is non-special unless either
\begin{itemize}
\item[(a)]
$k \geq 2$ and $\L'$ is nonempty and non-special,
or
\item[(b)] $\L'$ is special.
\end{itemize}
\end{corollary}

Finally one is able to make reductions
also in the case when $m_0 = d-m-1$.

\begin{proposition}
\label{m0=d-m-1}
Let $\L=\LS{d}{d-m-1}{b}{m}$ with $2 \leq m \leq d-1$.
Write $d = q(m-1)+\mu$ with $0\leq \mu \leq m-2$,
and $b = 2h+\epsilon$, with $\epsilon \in\{0,1\}$.
Then the system $\L$ is non-special of virtual dimension 
$d(m+2)-(b+1)m(m+1)/2$
unless
\begin{itemize}
\item[(a)] $q=h+1$, $\mu=\epsilon=0$, and $(m-1)(m+2) \geq 4h$,
in which case
\[
\dim \L = (m-1)(m+2)/2 - 2h,
\]
or
\item[(b)] $q = h$, $\epsilon = 0$, and $4q \leq \mu(\mu+3)$,
in which case $\dim \L = \mu(\mu+3)/2 - 2q$.
\end{itemize}
\end{proposition}

These statements above are taken directly from \cite{cm1};
for the argument below we need to extract the following
specific information.

\begin{corollary}
\label{highm0cor}
Let $\L=\LS{d}{d-m+k}{b}{m}$.
Suppose that $2 \leq m \leq d$ and $b=2h+1$ is odd.
If $-1 \leq k \leq 1$,
then $\L$ is non-special.
\end{corollary}

\section{The Induction Step for Large $d$}

Suppose now we want to prove the Main Conjecture
for systems $\LSH{d}{n}{m}$.
As noted above, the degeneration method
gives a bound (namely the dimensions of the limit system $\lo$)
for $\dim(\L)$
in terms of the four dimensions $\lP$, $\lPhat$, $\lF$, and $\lFhat$.
The last two are obtained with the results of the previous section.
The first two will be assumed to be non-special by induction.
With this approach we are able to prove the following theorem,
whose statement requires a bit of notation.

Define the function
\[
\dlow(\gamma,h) = 
\frac{\binom{m}{2}+\binom{\gamma+1}{2} + (2h+1)\binom{m+1}{2} -m\gamma -1}
{m+1-\gamma};
\]
notice that it is also a function of $m$, which we suppress.
Set
\[
D(m) = \max\{
\lfloor \frac{23m+16}{6}\rfloor,
\lceil \dlow(-1,\lceil \frac{m^2-1}{3m+4} \rceil) \rceil
\}.
\]
For $m$ large $D(m)$ is asymptotically $m^2/3$.

\begin{theorem}
\label{highd}
Fix $m \geq 2$ and let $D = D(m)$ as defined above.
Suppose that the Main Conjecture holds for all linear systems
$\LSH{d}{n}{m}$ with $d < D$.
Then the Main Conjecture holds for all linear systems
$\LSH{d}{n}{m}$.
\end{theorem}

\begin{pf}
We go by induction on $d$,
and therefore fix a $d \geq D(m)$,
and we assume that the Main Conjecture holds for all
linear systems $\LSH{d'}{n}{m}$ with $d' < d$.
We first take up the case when $\v < 0$,
and must show then that the system $\LSH{d}{n}{m}$ is empty,
because when $d \geq D(m)$ there are no $(-1)$-special systems.

To prove that this system is empty,
it suffices to find an $(a,b)$-degeneration
such that the limit dimension $\lo = -1$.
By Theorem \ref{dimLo}(a),
it suffices to have
$\lP + \lF \leq a-1$ and $\lPhat = \lFhat = -1$.

We take care to choose $b=2h+1$ odd,
and $a \in \{d-m-1,d-m,d-m+1\}$,
so that by Corollary \ref{highm0cor}
the system $\LF = \LS{d}{a}{b}{m}$ is non-special.

Since $d \geq D(m)$, in particular we have
$a-1 > (17m-2)/6$, so that by Theorem \ref{list-1special}
the systems $\LP$ and $\LPhat$ are not $(-1)$-special,
and therefore by induction they are not special.

We first claim that with these assumptions,
if $\lPhat = \lFhat = -1$,
the condition that $\lP + \lF \leq a-1$
(which is equivalent to $\rP+\rF \leq a-1$)
is automatic.
Indeed, if either one of the systems $\LP$ or $\LF$ are empty,
then it is obvious (because the restricted systems
whose dimensions are $\rP$ and $\rF$ 
both are subsystems of the complete system $|\O_R(a)|$
on the double curve $R$).
If neither is empty, then because they are non-special,
we have $\lP = \vP$ and $\lF = \vF$ so that
by Lemma \ref{vP+vF} we see that
$\lP+\lF = v + a \leq a - 1$ since $\v \leq -1$ in this case.

We are left with imposing that $\lPhat = \lFhat = -1$.
Let us write $a = d-m+\gamma$ with $\gamma \in \{-1,0,1\}$.
For the system $\LFhat = \LS{d}{a+1}{b}{m}$,
we remark that the $b$ lines through $p_0$ and the $p_i$'s
each split off $\gamma+1$ times,
leaving the residual system
$\LS{d-b(\gamma+1)}{d-m+\gamma+1-b(\gamma+1)}{b}{m-(\gamma+1)}$,
which by Corollary \ref{highm0cor} is non-special since $b$ is odd.
Therefore $\lFhat = -1$ if the virtual dimension
of this residual system is negative.
This lead to the inequality
\[
d \leq \dhigh(\gamma,h) = m +h m - 1 + h + h \gamma.
\]

For the system $\LPhat = \LSH{a-1}{n-b}{m}$,
we simply impose that $\vPhat \leq \v$, which gives the inequality
\[
d \geq \dlow(\gamma,h).
\]
We therefore obtain an inductive proof for this $d$
if we are able to choose $\gamma$ and $b=2h+1$
with $d$ in the interval $[\dlow(\gamma,h),\dhigh(\gamma,h)]$.

Both $\dlow$ and $\dhigh$ are increasing functions of $\gamma$.
It is a remarkable fact that
\[
\dhigh(-1,h) = \dlow(0,h) \;\;\text{ and }\;\; \dhigh(0,h) = \dlow(1,h)
\]
for every $h$.
Therefore for a fixed $h$ the three intervals
(given by the three values of $\gamma$)
match perfectly to give a single larger interval,
and if $d$ is in this larger interval,
we have executed the inductive proof.
This larger interval is $[\dlow(-1,h),\dhigh(1,h)]$.

We now want to vary $h$, and show that all integral $d$'s
larger than $D(m)$ are in one of these larger intervals.
This will be the case if the gap between the upper end
of one interval $\dhigh(1,h)$
is within one of the lower end of the next interval
$\dlow(-1,h+1)$,
since $\dhigh(1,h)$ is an integer.
This condition is that the difference is at most one,
i.e.
\begin{align*}
1 &\geq \dlow(-1,h+1) - \dhigh(1,h) \\
& = \frac{m^2+m+1-3mh-4h}{m+2}, \\
\end{align*}
which is equivalent to
\begin{equation}
\label{hinequality}
h \geq \frac{m^2-1}{3m+4}.
\end{equation}
We note in passing that $h$ \emph{can} be chosen to be this big;
this requires that $n$ be large enough, which it is since
$\v < 0$.  Specifically, we need $b=2h+1 < n$,
so that it is enough if $n \geq (2m^2-2)/(3m+4)$;
since $d(d+3)<nm(m+1)$, and $d \geq D(m)$, this is guaranteed.

Therefore as soon as $h$ is this big, there are no
integers $d$ which fail to be in one of the desired intervals.
These intervals begin at
\[
\lceil \dlow(-1,\lceil \frac{m^2-1}{3m+4} \rceil) \rceil,
\]
and since $D(m)$ is at least this quantity by assumption,
we are done.

We must now address the case
when the virtual dimension $\v$ is non-negative,
and we must show that the actual dimension $\l$
is equal to the expected dimension $\v$.
For this we assume that $\v$ is non-negative,
but that for this $d$ and $m$ the $n$ is maximal with $\v \geq 0$.
If we are able to prove that for this $d$, $n$, and $m$ we have
$\v = \l$, then for all smaller $n$'s we will also have $\v=\l$:
if the conditions imposed by $n$ general multiple points are
independent, then the conditions imposed by any fewer points are.

It suffices to find an $(a,b)$-degeneration
such that the limit dimension $\lo = \v$.
By Theorem \ref{dimLo}(b),
it suffices to have
$\rP + \rF \geq a-1$ and $\lP + \lF - a = \v$.

We will try to find $a$ and $b$
such that $\LP$ and $\LF$ are non-special
with virtual dimension at least $-1$;
if this is the case, then
$\lP = \vP$ and $\lF = \vF$,
and so we will obtain $\lP + \lF - a = \v$
automatically using Lemma \ref{vP+vF}.

Again we will take $b=2h+1$ odd,
and $a = d-m+\gamma$ with $\gamma \in \{-1,0,1\}$,
so that by Corollary \ref{highm0cor}
the system $\LF = \LS{d}{a}{b}{m}$ is non-special.
Its virtual dimension is
\[
\vF = d(m+1-\gamma) - \binom{\gamma+1}{2} - m^2+m\gamma-hm^2-hm
\]
which we need to be at least $-1$.
Requiring $\vF \geq -1$ is, by Lemma \ref{vP+vF},
equivalent to requiring $\vPhat \leq v$,
which we have already noted above in the first part of the proof
is the inequality
\[
d \geq \dlow(\gamma,h).
\]

We now also require that $\LFhat$ is empty;
the same analysis as we did in the $\v < 0$ case above
shows that this is implied by
\[
d \leq \dhigh(\gamma,h) = m +h m - 1 + h + h \gamma.
\]
Since $\lFhat \geq \vFhat$,
this will also imply that $\vFhat < 0$;
by Lemma \ref{vP+vF}, this also imposes that $\vP \geq \v$.

Since $d \geq D(m)$, in particular we have
$a-1 > (17m-2)/6$, so that by Theorem \ref{list-1special}
the systems $\LP$ and $\LPhat$ are not $(-1)$-special,
and therefore by induction they are not special.
In particular because $\vP \geq \v \geq 0$,
we have that $\LP$ is non-empty and non-special.

At this point we have $\lF + \lP -a = \v$,
so we need only address the inequality $\rP+\rF\geq a-1$,
which is equivalent to $\v \geq \lPhat + \lFhat + 1 = \lPhat$.

If $\lPhat = -1$, this inequality holds.
Otherwise since $\LPhat$ is non-special,
$\lPhat = \vPhat$, and the inequality follows from Lemma \ref{vP+vF}.

We therefore obtain an inductive proof for this $d$
if we are able to choose $\gamma$ and $b$
with $d$ in the interval $[\dlow(\gamma,h),\dhigh(\gamma,h)]$.
At this point the reader sees that the proof goes
identically as in the $\v < 0$ case,
and we are finished.

The only point to check is that we can take $h$ to be large enough
to satisfy (\ref{hinequality}).
Again this involves an inequality on $n$,
and here (and only here) we use the assumption that we have the
maximum $n$ with $\v \geq 0$.  The reader can check that this
is sufficient.
\end{pf}

By analysing the above proof, one sees that in the induction,
we only use that $\LP$ and $\LPhat$ are non-special,
and for this we need that the Main Conjecture holds
for the systems $\LSH{d-m-2}{n-b}{m}$, $\LSH{d-m-1}{n-b}{m}$,
$\LSH{d-m}{n-b}{m}$, and $\LSH{d-m+1}{n-b}{m}$.
This remark allows us to deduce the following.

\begin{corollary}
Fix $m \geq 2$ and let $D = D(m)$ as defined above.
Suppose that there is an $N \geq D(m)-1$
such that the Main Conjecture holds for linear systems
$\LSH{d}{n}{m}$ with $N-m-1 \leq d \leq N$.
Then the Main Conjecture holds for all linear systems
$\LSH{d}{n}{m}$ with $d \geq N$.
\end{corollary}

\section{Proof of the Main Conjecture for $m \leq 12$}
Suppose we want to prove the main conjecture
for a fixed degree $d$ and fixed multiplicity $m$
and all numbers of points $n$.
There is a critical number $n_0 = n_0(d,m)$
such that the virtual dimension of $\LSH{d}{n_0}{m}$
is positive, but that of $\LSH{d}{n_0+1}{m}$
is negative.
If one can show that the critical system
$\LSH{d}{n_0}{m}$ is non-special,
and that the system $\LSH{d}{n_0+1}{m}$ is empty,
then $\LSH{d}{n}{m}$ will be non-special for all $n$.

\begin{proposition}
Fix $m$ and $d \leq 3m$.
Then for all $n$ the Main Conjecture holds
for the homogeneous linear system $\LSH{d}{n}{m}$.
\end{proposition}

\begin{pf}
We know that the Main Conjecture holds
for all $n \leq 9$.
Therefore if $d$ and $m$ are fixed with $n_0(d,m) \leq 8$,
then the Main Conjecture will hold for all $n$.
This is true if $d < 3m$.

If $d = 3m$, then $n_0 = 9$, and the linear system
$\LSH{3m}{9}{m}$ is non-special, of dimension $0$,
consisting of the unique multiple cubic through the
$9$ general points.
Therefore the system $\LSH{3m}{10}{m}$ is empty.
\end{pf}

Using Theorem \ref{highd} and the previous Proposition,
the Main Conjecture for a fixed $m$
and all $d$ and $n$ will follow
if one can show that for all $d$ in the interval
$[3m+1,D(m)-1]$,
the system $\LSH{d}{n_0}{m}$ is non-special
and the system $\LSH{d}{n_0+1}{m}$
is empty.

These intervals for $m \leq 12$ are given in the table below.
\begin{center}
\begin{tabular}{r|ccccccccccc}
$m$      & 2 &  3 &  4 &  5 &  6 &  7 &  8 &  9 & 10 & 11 & 12 \\ \hline
$3m+1$   & 7 & 10 & 13 & 16 & 19 & 22 & 25 & 28 & 31 & 34 & 37 \\ \hline
$D(m)-1$ & 9 & 13 & 17 & 20 & 24 & 28 & 32 & 36 & 40 & 50 & 55
\end{tabular}
\end{center}

We have shown in Theorem \ref{highd} that the degeneration
method will always work for investigating the system
$\LSH{d}{n}{m}$ when $d \geq D(m)$,
reducing the computation to the knowledge of the dimensions
of homogeneous systems with lower $d$ and $n$.
For smaller values of $d$ the method often works anyway;
the estimates given in the proof of the Theorem
simply do not guarantee a suitable $(a,b)$-degeneration,
but usually one exists, even for $d$'s in the middle
range $[3m+1,D(m)-1]$.

We have written a computer program to investigate,
for a fixed $m$,
all $d$'s in this middle range,
and the two critical values of $n$, namely
$n_0$ and $n_0+1$,
searching for a suitable $a$ and $b$ to execute
the recursion which comes out of the $(a,b)$-degeneration.
If the program successfully finds an $a$ and $b$,
it outputs the values and goes on to the next case.
If the program does not find any suitable $a$ and $b$,
the value of $d$ and the critical $n$ for which the recursion fails
is printed.

Below we present in tabular form, for each $m$ with $2 \leq m \leq 12$,
for each $d$ with $3m+1 \leq d \leq D(m)-1$,
and for the critical values $n_0$ and $n_0+1$,
a suitable $a$ and $b$ for which the recursion succeeds.
We leave it to the reader to check the details of these finitely many
computations.
We note that when the virtual dimension of $\LSH{d}{n_0}{m}$
is exactly $-1$,
a proof for that system suffices for the $n_0+1$ system as well.

\[
\begin{array}{ccc}
\begin{tabular}{ccccc}
$d$ & $n$ & $m$ & $v$ & $(a,b)$ \\ \hline
7 & 12 & 2 & -1 & (4,8) \\
8 & 15 & 2 & -1 & (5,9) \\
9 & 18 & 2 & 0 & (6,9) \\
9 & 19 & 2 & -3 & (6,12) \\
10 & 11 & 3 & -1 & (6,7) \\
11 & 13 & 3 & -1 & (7,7) \\
12 & 15 & 3 & 0 & (8,8) \\
13 & 16 & 3 & -6 & (9,8)\\
13 & 10 & 4 & 4 & (7,7)\\
13 & 11 & 4 & -6 & (7,7)\\
14 & 12 & 4 & -1 & (8,8)\\
15 & 13 & 4 & 5 & (9,9)\\
15 & 14 & 4 & -5 & (10,7)\\
16 & 15 & 4 & 2 & (11,8)\\
16 & 16 & 4 & -8 & (11,9)\\
17 & 17 & 4 & 0 & (11,10)\\
17 & 18 & 4 & -10 & (11,11)
\end{tabular}
&
\begin{tabular}{ccccc}
$d$ & $n$ & $m$ & $v$ & $(a,b)$ \\ \hline
16 & 10 & 5 & 2 & (9,7) \\
16 & 11 & 5 & -13 & (9,7)\\
17 & 11 & 5 & 5 & (11,7)\\
17 & 12 & 5 & -10 & (10,8)\\
18 & 12 & 5 & 9 & (11,8)\\
18 & 13 & 5 & -6 & (12,7)\\
19 & 14 & 5 & -1 & (14,7)\\
20 & 15 & 5 & 5 & (14,7)\\
20 & 16 & 5 & -10 & (14,9) \\
19 & 10 & 6 & -1 & FAIL\\
20 & 11 & 6 & -1 & (12,7)\\
21 & 12 & 6 & 0 & (16,5)\\
21 & 13 & 6 & -21 & (12,9)\\
22 & 13 & 6 & 2 & (15,7)\\
22 & 14 & 6 & -19 & (15,7)\\
23 & 14 & 6 & 5 & (16,7)\\
23 & 15 & 6 & -16 & (16,8)
\end{tabular}
&
\begin{tabular}{ccccc}
$d$ & $n$ & $m$ & $v$ & $(a,b)$ \\ \hline
24 & 15 & 6 & 9 & (16,8)\\
24 & 16 & 6 & -12 & (16,9)\\
22 & 9 & 7 & 23 & (16,5)\\
22 & 10 & 7 & -5 & FAIL\\
23 & 10 & 7 & 19 & (17,5)\\
23 & 11 & 7 & -9 & (17,5)\\
24 & 11 & 7 & 16 & (18,5)\\
24 & 12 & 7 & -12 & (19,5)\\
25 & 12 & 7 & 14 & (19,5)\\
25 & 13 & 7 & -14 & (19,6)\\
26 & 13 & 7 & 13 & (18,7)\\
26 & 14 & 7 & -15 & (18,7)\\
27 & 14 & 7 & 13 & (19,7)\\
27 & 15 & 7 & -15 & 19,7)\\
28 & 15 & 7 & 14 & (19,8)\\
28 & 16 & 7 & -14 & (19,9)\\
   &    &   &     &
\end{tabular}
\end{array}
\]

\[
\begin{array}{ccc}
\begin{tabular}{ccccc}
$d$ & $n$ & $m$ & $v$ & $(a,b)$ \\ \hline
25 & 9 & 8 & 26 & (18,5)\\
25 & 10 & 8 & -10 & FAIL\\
26 & 10 & 8 & 17 & (15,7)\\
26 & 11 & 8 & -19 & (19,5)\\
27 & 11 & 8 & 9 & (20,5)\\
27 & 12 & 8 & -27 & (16,8)\\
28 & 12 & 8 & 2 & (16,8)\\
28 & 13 & 8 & -34 & (15,9)\\
29 & 12 & 8 & 32 & (20,6)\\
29 & 13 & 8 & -4 & (20,7)\\
30 & 13 & 8 & 27 & (21,6)\\
30 & 14 & 8 & -9 & (21,7)\\
31 & 14 & 8 & 23 & (22,7)\\
31 & 15 & 8 & -13 & (20,9)\\
32 & 15 & 8 & 20 & (21,9)\\
32 & 16 & 8 & -16 & (21,9)\\
28 & 9 & 9 & 29 & (20,5)\\
28 & 10 & 9 & -19 & FAIL
\end{tabular}
&
\begin{tabular}{ccccc}
$d$ & $n$ & $m$ & $v$ & $(a,b)$ \\ \hline
29 & 10 & 9 & 14 & FAIL\\
29 & 11 & 9 & -31 & (21,5)\\
30 & 11 & 9 & 0 & (18,7)\\
30 & 12 & 9 & -45 & (17,8)\\
31 & 11 & 9 & 32 & (20,7)\\
31 & 12 & 9 & -13 & (18,8)\\
32 & 12 & 9 & 20 & (24,5)\\
32 & 13 & 9 & -25 & (18,9)\\
33 & 13 & 9 & 9 & (22,7)\\
33 & 14 & 9 & -36 & (22,8)\\
34 & 14 & 9 & -1 & (24,7)\\
35 & 14 & 9 & 35 & (26,6)\\
35 & 15 & 9 & -10 & (25,7)\\
36 & 15 & 9 & 27 & (24,8)\\
36 & 16 & 9 & -18 & (24,9)\\
31 & 9 & 10 & 32 & (21,5)\\
31 & 10 & 10 & -23 & FAIL\\
32 & 10 & 10 & 10 & FAIL
\end{tabular}
&
\begin{tabular}{ccccc}
$d$ & $n$ & $m$ & $v$ & $(a,b)$ \\ \hline
32 & 11 & 10 & -45 & (19,7)\\
33 & 10 & 10 & 44 & (19,7)\\
33 & 11 & 10 & -11 & (20,7)\\
34 & 11 & 10 & 24 & (21,7)\\
34 & 12 & 10 & -31 & (19,8)\\
35 & 12 & 10 & 5 & (20,8)\\
35 & 13 & 10 & -50 & (19,9)\\
36 & 12 & 10 & 42 & (27,5)\\
36 & 13 & 10 & -13 & (24,7)\\
37 & 13 & 10 & 25 & (25,7)\\
37 & 14 & 10 & -30 & (26,7)\\
38 & 14 & 10 & 9 & (27,7)\\
38 & 15 & 10 & -46 & (24,9)\\
39 & 14 & 10 & 49 & (28,7)\\
39 & 15 & 10 & -6 & (28,7)\\
40 & 15 & 10 & 35 & (27,8)\\
40 & 16 & 10 & -20 & (27,9)\\
   &    &   &     &
\end{tabular}
\end{array}
\]

\[
\begin{array}{ccc}
\begin{tabular}{ccccc}
$d$ & $n$ & $m$ & $v$ & $(a,b)$ \\ \hline
34 & 9 & 11 & 35 & (23,5)\\
34 & 10 & 11 & -31 & FAIL\\
35 & 10 & 11 & 5 & FAIL\\
35 & 11 & 11 & -61 & (22,7)\\
36 & 10 & 11 & 42 & (21,7)\\
36 & 11 & 11 & -24 & (22,7)\\
37 & 11 & 11 & 14 & (23,7)\\
37 & 12 & 11 & -52 & (27,6)\\
38 & 11 & 11 & 53 & (27,5)\\
38 & 12 & 11 & -13 & (22,8)\\
39 & 12 & 11 & 27 & (23,8)\\
39 & 13 & 11 & -39 & (26,7)\\
40 & 13 & 11 & 2 & (27,7)\\
40 & 14 & 11 & -64 & (26,8)\\
41 & 13 & 11 & 44 & (28,7)\\
41 & 14 & 11 & -22 & (29,7)\\
42 & 14 & 11 & 21 & (27,8)\\
42 & 15 & 11 & -45 & (27,9)\\
43 & 15 & 11 & -1 & (31,7)\\
44 & 15 & 11 & 44 & (29,9)\\
44 & 16 & 11 & -22 & (29,9)\\
45 & 16 & 11 & 24 & (30,9)\\
45 & 17 & 11 & -42 & (31,9)
\end{tabular}
&
\begin{tabular}{ccccc}
$d$ & $n$ & $m$ & $v$ & $(a,b)$ \\ \hline
46 & 17 & 11 & 5 & (32,9)\\
46 & 18 & 11 & -61 & (29,11)\\
47 & 17 & 11 & 53 & (31,10)\\
47 & 18 & 11 & -13 & (32,10)\\
48 & 18 & 11 & 36 & (35,9)\\
48 & 19 & 11 & -30 & (32,11)\\
49 & 19 & 11 & 20 & (34,10)\\
49 & 20 & 11 & -46 & (35,10)\\
50 & 20 & 11 & 5 & (35,10)\\
50 & 21 & 11 & -61 & (35,11)\\
37 & 9 & 12 & 38 & (25,5)\\
37 & 10 & 12 & -40 & FAIL\\
38 & 10 & 12 & -1 & FAIL\\
39 & 10 & 12 & 39 & FAIL\\
39 & 11 & 12 & -39 & (28,5)\\
40 & 11 & 12 & 2 & (24,7)\\
40 & 12 & 12 & -76 & (29,6)\\
41 & 11 & 12 & 44 & (25,7)\\
41 & 12 & 12 & -34 & (24,8)\\
42 & 12 & 12 & 9 & (24,8)\\
42 & 13 & 12 & -69 & (23,9)\\
43 & 12 & 12 & 53 & (25,8)\\
43 & 13 & 12 & -25 & (29,7)
\end{tabular}
&
\begin{tabular}{ccccc}
$d$ & $n$ & $m$ & $v$ & $(a,b)$ \\ \hline
44 & 13 & 12 & 20 & (30,7)\\
44 & 14 & 12 & -58 & (28,8)\\
45 & 13 & 12 & 66 & (31,7)\\
45 & 14 & 12 & -12 & (32,7)\\
46 & 14 & 12 & 35 & (33,7)\\
46 & 15 & 12 & -43 & (29,9)\\
47 & 15 & 12 & 5 & (34,7)\\
47 & 16 & 12 & -73 & (31,9)\\
48 & 15 & 12 & 54 & (32,8)\\
48 & 16 & 12 & -24 & (32,9)\\
49 & 16 & 12 & 26 & (33,9)\\
49 & 17 & 12 & -52 & (32,10)\\
50 & 17 & 12 & -1 & (34,9)\\
51 & 17 & 12 & 51 & (36,9)\\
51 & 18 & 12 & -27 & (34,10)\\
52 & 18 & 12 & 26 & (35,10)\\
52 & 19 & 12 & -52 & (35,10)\\
53 & 19 & 12 & 2 & (39,9)\\
53 & 20 & 12 & -76 & (38,10)\\
54 & 19 & 12 & 57 & (40,8)\\
54 & 20 & 12 & -21 & (38,10)\\
55 & 20 & 12 & 35 & (39,10)\\
55 & 21 & 12 & -43 & (39,11)
\end{tabular}
\end{array}
\]

For those linear systems for which the method fails,
one must argue in a different way;
if one can successfully show that
these finitely many systems are non-special,
then one has proved that all homogeneous linear systems
$\LSH{d}{n}{m}$ with this fixed multiplicity $m$
satisfy the Main Conjecture.
These finitely many systems for $m \leq 12$ which we
must deal with are presented below.

\begin{center}
\begin{tabular}{ccccl}
$d$ & $n$ & $m$ & $\v$ & reason for truth of the Main Conjecture \\ \hline
19 & 10 & 6 & -1 & \cite{hirschowitz1}, or implied by $\LSH{38}{10}{12}$ \\ 
22 & 10 & 7 & -5 & restrict $\LSH{38}{10}{12}$ \\ 
25 & 10 & 8 & -10 & $D-3K$, or restrict $\LSH{38}{10}{12}$ \\ 
28 & 10 & 9 & -16 & $D-4K$, or restrict $\LSH{38}{10}{12}$ \\ 
29 & 10 & 9 & 14 & $D-4K+H$ \\ 
31 & 10 & 10 & -23 & $D-5K$, or restrict $\LSH{38}{10}{12}$ \\ 
32 & 10 & 10 & 10 & $D-5K+H$ \\ 
34 & 10 & 11 & -31 & $D-6K$, or restrict $\LSH{38}{10}{12}$ \\ 
35 & 10 & 11 & 5 & $D-6K+H$ \\ 
37 & 10 & 12 & -40 & $D-7K$, or implied by (or restrict) $\LSH{38}{10}{12}$ \\ 
38 & 10 & 12 & -1 & Gimigliano's Thesis \cite{gimigliano1} \\ 
39 & 10 & 12 & 39 & implied by $\LSH{38}{10}{12}$ \\
\end{tabular}
\end{center}

In the last column of the table above,
for the reader's convenience,
we have given an abbreviated description
of the argument or the reference used below
to show that these systems all satisfy the Main Conjecture
(and hence all have the expected dimension).
Notationally, $D$ represents
a general element of the system $\LSH{16}{10}{5}$,
$K$ is the canonical class, and $H$ is the line class.

\begin{theorem}
For every linear system $\LSH{d}{n}{m}$
with $m \leq 12$,
the Main Conjecture is true.
\end{theorem}

\begin{pf}
For $m\leq 5$, the degeneration method
works in every case, and there is nothing more to do.
For $6 \leq m \leq 12$, the degeneration method
reduces us to checking each of the $12$ cases
presented above.
These we take up in turn.

Firstly, in the thesis of Gimigliano \cite{gimigliano1},
he uses the Horace Method developed in
\cite{hirschowitz1} to prove that the system
$\L = \LSH{38}{10}{12}$ is empty.
Since its virtual dimension is $-1$, this implies
that both $H^0$ and $H^1$ of this system
(as a complete linear system on the blowup of the plane
at the $10$ general points)
are zero.
If we denote by $H$ the line class on blowup,
then also of course the system $\L-H = \LSH{37}{10}{12}$
must also be empty.

Consider the system $\L+H = \LSH{39}{10}{12}$.
Restricting this to a general line gives the short exact sequence
\[
0 \to \L \to \L+H \to \O_H(39) \to 0,
\]
and since we have that $H^1(\L) = H^1(\O_H(39)) = 0$,
we see that also $H^1(\L+H) = 0$;
therefore $\L+H$ is non-special, of dimension $39$
as expected.

The system $\LSH{19}{10}{6}$ must also be empty,
since if it contained an effective divisor $F$,
then $2F$ would be a member of $\LSH{39}{10}{12}$.
This system $\LSH{19}{10}{6}$ was also shown to be empty
by the Horace Method in \cite{hirschowitz1}.
We note for the argument below that since the virtual dimension
of $\LSH{19}{10}{6}$ is $-1$,
its emptyness implies that its $H^1$ is zero.

Finally the emptyness of $\LSH{38}{10}{12}$
also implies the emptyness of several other systems
of the form $\LSH{d}{10}{m}$ with negative virtual dimension $\v$.
Suppose on the contrary that an effective divisor $C$
existed in the system $\LSH{d}{10}{m}$.
Taking the ideal sequence of $C$
and twisting by the line bundle corresponding
to the linear system $\LSH{38}{10}{12}$ gives
\[
0 \to \LSH{38-d}{10}{(12-m)} \to \LSH{38}{10}{12} \to \LSH{38}{10}{12}|_C \to 0
\]
which gives a contradiction if the virtual dimension
of the kernel system $\LSH{38-d}{10}{(12-m)}$ is non-negative;
then the sheaf on the left will have sections, while that
in the middle does not.
This argument is successful for proving the emptyness
of $\LSH{22}{10}{7}$ (where the kernel system is $\LSH{16}{10}{5}$),
$\LSH{25}{10}{8}$ (the kernel system is $\LSH{13}{10}{4}$),
$\LSH{28}{10}{9}$ (the kernel system is $\LSH{10}{10}{3}$),
$\LSH{31}{10}{10}$ (the kernel system is $\LSH{7}{10}{2}$),
and $\LSH{34}{10}{11}$ (the kernel system is $\LSH{4}{10}{1}$).

An alternate argument for most of these empty systems
is as follows.
Consider the system $\LSH{16}{10}{5}$,
which is non-special of dimension $2$;
let $D$ be a general divisor in this system.
Note that $D$ is irreducible,
since if not, by symmetry, any irreducible component
would have to generate (under the permutation group of the $10$ points)
a homogeneous linear system $\LSH{\delta}{10}{\mu}$,
all of whose members were components of $D$.
Since we would have $\mu \leq 5$,
and these systems would have to have non-negative expected dimension,
we would need
$\delta \geq 4$ if $\mu = 1$,
$\delta \geq 7$ if $\mu = 2$,
$\delta \geq 10$ if $\mu = 3$, and
$\delta \geq 13$ if $\mu = 4$.
But then $\dim \LSH{\delta}{10}{\mu}$ would be at least $3$,
contradicting the fact that the dimension of $|D|$ is $2$.

Consider the systems $D-mK$;
since $(D-mK)\cdot D = 6-2m$, and $D$ is irreducible and effective
with $D^2 > 0$ and $\dim|D| = 2$,
we see that $|D-mK|$ is empty as soon as $m \geq 3$.
This shows that the systems
$\LSH{25}{10}{8}$ (m=3),
$\LSH{28}{10}{9}$ (m=4),
$\LSH{31}{10}{10}$ (m=5),
$\LSH{34}{10}{11}$ (m=6), and
$\LSH{37}{10}{12}$ (m=7)
are all empty.

The above arguments handle all of the cases
where the virtual dimension is negative.
We now turn to the remaining three cases
where the virtual dimension is positive.

Consider $H - (m-1)K = \LSH{3m-2}{10}{(m-1)}$;
its virtual dimension is non-negative for $m \leq 6$.
Let $C$ be a general member of this system.
An argument as above shows that $C$ is irreducible;
denote by $g$ its arithmetic genus.
Specifically, these systems are

\begin{center}
\begin{tabular}{cccc}
$m$ & $H-(m-1)K$ & $g$ & $(D-mK+H)\cdot C$ \\ \hline
$3$ & $\LSH{7}{10}{2}$ & $5$ & $22$ \\
$4$ & $\LSH{10}{10}{3}$ & $6$ & $20$ \\
$5$ & $\LSH{13}{10}{4}$ & $6$ & $16$ \\
$6$ & $\LSH{16}{10}{5}$ & $5$ & $10$
\end{tabular}
\end{center}

Note that in all cases $(D-mK+H)\cdot C$ is at least $2g$,
so that the restricted system $(D-mK+H)|_C$ is non-special on $C$.
Using the exact sequence
\[
0 \to D-K = \LSH{19}{10}{6} \to D-mK+H \to (D-mK+H)|_C \to 0
\]
we see that since also $H^1(\LSH{19}{10}{6})=0$,
we have $H^1(D-mK+H) = 0$ too,
proving that these four systems are non-special.
For $3 \leq m \leq 6$ these systems are
$\LSH{26}{10}{8}$ (m=3),
$\LSH{29}{10}{9}$ (m=4),
$\LSH{32}{10}{10}$ (m=5), and
$\LSH{35}{10}{11}$ (m=6);
this provides an alternate proof for $\LSH{26}{10}{8}$
for which the degeneration method also worked
(in fact using $a=15$ and $b=7$).

This completes the analysis of all systems
for which the degeneration method failed,
and finishes the proof.
\end{pf}

\end{document}